# Simple and fast Algorithm for Finding Roots of Error-Locator Polynomials: Modulus Search.

**Gennady N. Glushchenko[1]**

[1]Moscow Institute of Physics and Technology, 9 Institutskiy per., Dolgoprudny, Moscow Region, 141701, Russia

e-mail: glushchenko.gn@mipt.ru

**ABSTRACT** A novel very simple method for finding roots of polynomials over finite fields has been proposed. The essence of the proposed method is to search the roots via nested cycles over the subgroups of the multiplicative group of the Galois field. The modified Chien search is actually used in the inner cycles, but the internal polynomials are small. The word "modulus" was used because the search is doing on subsets like $\alpha^{a+bi}$, where a,b=const. In addition, modulo division of polynomials is actively used. The algorithm is applicable not for all Galois fields, but for selective ones, starting from $GF(2^8)$. The algorithm has an advantage for large polynomials. The number of operations is significant for small polynomials, but it grows very slowly with the degree of the polynomial. When the polynomial is large or very large, the proposed method can be 10-100 times faster than Chien search.

**INDEX TERMS** BCH codes, decoding, error correction codes, error locator polynomial, Galois fields, Reed-Solomon codes.

## I. INTRODUCTION

It is well known that one of the most time-consuming stages of decoding process of Reed-Solomon, BCH and some other codes is finding roots of the error-locator polynomial. The most widely known root finding algorithm is Chien search method. It has very high time complexity for the case of large fields and polynomials of high degree. There are many algorithms that speed up the root search for small fields. I propose a very simple algorithm resembling redesigned Chien search that can be very effective for large fields and polynomials.

## II. ABBREVIATIONS

We work with Galois field $GF(2^n)$.
The number of nonzero elements is $m = 2^n - 1$.
Our polynomial is p,

$$p(x) = \sum_{i=0}^{t} p_i x^i$$

It's degree $t = \deg(p)$.
$\alpha$ is a primitive element of $GF(2^n)$.

## III. CHIEN SEARCH

We will start with a little investigation of Chien search. The most straightforward form of the root search is:

**Algorithm 1:** rough search

given polynomial p
for i = 0 to m - 1 do
  if $p(\alpha^i)$ == 0 then
    output($\alpha^i$)
  end if
end for

Let's introduce operation

$$p(x) := p(\alpha \cdot x)$$

which is a synonym to

$$p_i := \alpha^i p_i$$

for all i from 1 to deg(p).
We can write Chien search algorithm in the following form:

**Algorithm 2:** Chien search

given polynomial p
for i = 0 to m - 1 do
  if i > 0 then
    $p(x) := p(\alpha \cdot x)$
  end if
  if p(1) == 0 then
    output($\alpha^i$)
  end if
end for

The finding of p(1) does not require any multiplications and is performed very quickly. Total complexity of the last algorithm is (m-1)·t multiplications.

Finally, we can add the following optimization:

**Algorithm 3:** Chien search with corner exclusion

```
for i = 0 to m - 1 do
  if i > 0 then
    p(x) := p(α·x)
  end if
  if p(1) == 0 then
    output(α^i)
    p := p / (x + 1)    // optimization
  end if
end for
```

The essence of this optimization is in excluding each obtained root from the polynomial. The polynomial becomes smaller and the further cycles are performed faster. By itself, the operation p := p / (x + 1) is done quickly. However, the operation p(x) := p(α·x) usually assumes constant coefficients of the polynomial. We can prepare discrete logarithms of $p_i$ in advance, and then only take the corresponding exponentials. If the polynomial changes, then we need to restart the Chien search again. This gives an additional number of operations ≈$t^2$. The total number of operations for typical error-locator polynomials will change from ≈m·t to ≈m·t/2+$t^2$. For large values of t=deg(p), there is not much sense in such optimization. But we will use it further in the inner loop, where the polynomial has a small degree.

## IV. THE PROPOSED ALGORITHM

Suppose that we can split m into multipliers

$$m = m1·m2$$

We can write a loop variable i in the form i=i1+m1·i2 and rewrite the rough algorithm 1 in the form of two-cycle search:

**Algorithm 4:** two-cycle rough search

```
given polynomial p
for i1 = 0 to m1 - 1 do
  for i2 = 0 to m2 - 1 do
    if p(α^(i1+m1*i2)) == 0 then
      output(α^(i1+m1*i2))
    end if
  end for
end for
```

Note that the values $α^{i1+m1·i2}$ in the inner loop satisfy the equation

$$x^{m2} + c = 0,$$

where $c = α^{m2·i1}$ is a constant for the inner loop.

Therefore, we can replace p with p mod($x^{m2}$+c) in the inner loop:

**Algorithm 5:** two-cycle optimized search

```
given polynomial p
for i1 = 0 to m1 - 1 do
  p2 := p mod(x^m2 + α^(m2*i1))
  for i2 = 0 to m2 - 1 do
    if p2(α^(i1+m1*i2)) == 0 then
      output(α^(i1+m1*i2))
    end if
  end for
end for
```

We have introduced an intermediate polynomial p2. His degree is deg(p2) ≤ m2-1. If deg(p) >> m2, then deg(p2) << deg(p) and we will get a significant acceleration in the inner cycle.

Finally, we can rewrite the last algorithm in 'Chien' form:

**Algorithm 6:** the proposed algorithm base scheme

```
given polynomial p
for i1 = 0 to m1 - 1 do
  if i1 > 0 then
    p(x) := p(α·x)
  end if
  p2 := p mod(x^m2 + 1)
  for i2 = 0 to m2 - 1 do
    if i2 > 0 then
      p2(x) := p2(α^m1·x)
    end if
    if p2(1) == 0 then
      output(α^(i1+m1*i2))
    end if
  end for
end for
```

The complexity of this algorithm is

$$(m1-1)·t + m1·(m2-1)·min(m2-1,t)$$

multiplications. This is a piecewise linear function of t. It can be noted that this function will be smaller for the case of m1>m2 rather than m1<m2.

We have described the basic idea of the proposed method. If m can be decomposed into several multipliers, then there may be several cycles, respectively. The highest performance is achieved if the maximum possible number of cycles is used and if they are ordered so that the inner cycles have the smallest periods. Also, for a small m, it makes sense to put optimization (Algorithm 3) in the innermost loop. This can give us an additional acceleration of up to several percent.

## V. GF($2^8$)

Here is the full version of the proposed algorithm for GF($2^8$).
Its scheme: m = 255 = 17 × 5 × 3

**Algorithm 7:** the proposed algorithm for GF($2^8$)

```
given polynomial p
for i1 = 0 to 17-1 do
  if i1 > 0 then
    p(x) := p(α·x)
  end if
  p2 := p mod(x^15 + 1)
  for i2 = 0 to 5-1 do
    if i2 > 0 then
      p2(x) := p2(α^17·x)
    end if
    p3 := p2 mod(x^3 + 1)
    for i3 = 0 to 3-1 do
      if i3 > 0 then
        p3(x) := p3(α^85·x)
      end if
      if p3(1) == 0 then
        output(α^(i1+17*i2+85*i3))
        p3 := p3 / (x + 1)
      end if
    end for
  end for
end for
```

TABLE I
OUR METHOD PERFORMANCE FOR GF($2^8$)

| t | Chien | max muls | mean muls | speedup |
|---|---|---|---|---|
| 10 | 2540 | 1180 | 1170 | 2.17 |
| 20 | 5080 | 1612 | 1592 | 3.19 |
| 30 | 7620 | 1772 | 1742 | 4.37 |
| 40 | 10160 | 1932 | 1892 | 5.37 |
| 50 | 12700 | 2092 | 2041 | 6.22 |
| 60 | 15240 | 2252 | 2192 | 6.95 |
| 70 | 17780 | 2412 | 2341 | 7.60 |
| 80 | 20320 | 2572 | 2490 | 8.16 |
| 90 | 22860 | 2732 | 2638 | 8.66 |
| 100 | 25400 | 2892 | 2787 | 9.11 |
| 110 | 27940 | 3052 | 2936 | 9.52 |
| 120 | 30480 | 3212 | 3085 | 9.88 |
| 130 | 33020 | 3372 | 3231 | 10.22 |
| 140 | 35560 | 3532 | 3378 | 10.53 |
| ≥14 | 254·t | 1292+16·t | | |

The number of additions is approximately equal to the number of multiplications, so we give only the latter. The average number of multiplications for the proposed method was obtained through computer simulation. First, the pseudorandom corners $\{x_i\}$ from GF($2^8$) were generated, nonzero and different. Then, the polynomial

$$p(x) = \prod_i (x - x_i)$$

was calculated. The proposed method was performed and checked. The calculation of Galois multiplications was performed through a debugging counter in the program code. Each test was called many times ($10^3$) and the average value of the operations was calculated.

It can be seen that the proposed method begins to demonstrate good efficiency only at large values of t. We can try to compare it with the best alternative methods for GF($2^8$) [3]. In general, direct comparison is difficult. Alternative methods for GF($2^8$) have a very small number of multiplications, mostly additions. It is necessary to somehow compare Galois multiplications with additions. If we take the coefficient mul/add=4, then the proposed method will have comparable performance at t≈100-150. With a higher mul/add ratio, the proposed method will work worse.

## VI. GALOIS FIELDS FROM GF($2^9$) TO GF($2^{16}$)

The full advantage of the proposed method is revealed for large and very large polynomials and for "convenient" Galois fields. The main criterion is that m can be decomposed into small prime factors. We present performance data for all fields from GF($2^9$) to GF($2^{16}$).

The cycle configurations, from external cycles to internal:
GF($2^9$): 73 × 7
GF($2^{10}$): 31 × 11 × 3
GF($2^{11}$): 89 × 23
GF($2^{12}$): 13 × 7 × 5 × 3 × 3
GF($2^{13}$): the method is not applicable
GF($2^{14}$): 127 × 43 × 3
GF($2^{15}$): 151 × 31 × 7
GF($2^{16}$): 257 × 17 × 5 × 3

TABLE II
OUR METHOD PERFORMANCE FOR GF($2^9$)

| t | Chien | max muls | mean muls | speedup |
|---|---|---|---|---|
| 10 | 5100 | 3348 | 3319 | 1.54 |
| 20 | 10200 | 4068 | 4008 | 2.54 |
| 30 | 15300 | 4788 | 4699 | 3.26 |
| 40 | 20400 | 5508 | 5387 | 3.79 |
| 50 | 25500 | 6228 | 6078 | 4.20 |
| 70 | 35700 | 7668 | 7457 | 4.79 |
| 100 | 51000 | 9828 | 9528 | 5.35 |
| 150 | 76500 | 13428 | 12977 | 5.89 |
| 200 | 102000 | 17028 | 16428 | 6.21 |
| 300 | 153000 | 24228 | 23299 | 6.57 |
| ≥6 | 510·t | 2628+72·t | | |

TABLE III
OUR METHOD PERFORMANCE FOR GF($2^{10}$)

| t | Chien | max muls | mean muls | speedup |
|---|---|---|---|---|
| 10 | 10220 | 4764 | 4754 | 2.15 |
| 20 | 20440 | 8164 | 8144 | 2.51 |
| 30 | 30660 | 11564 | 11537 | 2.66 |
| 40 | 40880 | 12484 | 12442 | 3.29 |
| 50 | 51100 | 12784 | 12734 | 4.01 |
| 70 | 71540 | 13384 | 13314 | 5.37 |
| 100 | 102200 | 14284 | 14184 | 7.21 |
| 150 | 154300 | 15784 | 15632 | 9.81 |
| 200 | 204400 | 17284 | 17083 | 11.97 |
| 300 | 306600 | 20284 | 19976 | 15.35 |
| 400 | 408800 | 23284 | 22863 | 17.88 |
| 500 | 511000 | 26284 | 25745 | 19.85 |
| ≥32 | 1022·t | 11284+30·t | | |

TABLE IV
OUR METHOD PERFORMANCE FOR GF($2^{11}$)

| t | Chien | max muls | mean muls | speedup |
|---|---|---|---|---|
| 10 | 20460 | 20460 | 20352 | 1.01 |
| 20 | 40920 | 40920 | 40696 | 1.01 |
| 30 | 61380 | 45716 | 45388 | 1.35 |
| 40 | 81840 | 46596 | 46161 | 1.77 |
| 50 | 102300 | 47476 | 46923 | 2.18 |
| 70 | 143220 | 49236 | 48464 | 2.96 |
| 100 | 204600 | 51876 | 50775 | 4.03 |
| 150 | 306900 | 56276 | 54624 | 5.62 |
| 200 | 409200 | 60676 | 58479 | 7.00 |
| 300 | 613800 | 69476 | 66171 | 9.28 |
| 400 | 818400 | 78276 | 73881 | 11.08 |
| 500 | 1023000 | 87076 | 81574 | 12.54 |
| 700 | 1432200 | 104676 | 96969 | 14.77 |
| 1000 | 2046000 | 131076 | 120078 | 17.04 |
| ≥22 | 2046·t | 43076+88·t | | |

TABLE V
OUR METHOD PERFORMANCE FOR GF($2^{12}$)

| t | Chien | max muls | mean muls | speedup |
|---|---|---|---|---|
| 10 | 40940 | 17280 | 17272 | 2.37 |
| 20 | 81880 | 21820 | 21802 | 3.76 |
| 30 | 122820 | 26360 | 26334 | 4.66 |
| 40 | 163760 | 30900 | 30864 | 5.31 |
| 50 | 204700 | 33256 | 33197 | 6.17 |
| 70 | 286580 | 35056 | 34989 | 8.19 |
| 100 | 409400 | 37756 | 37655 | 10.87 |
| 150 | 614100 | 42256 | 42106 | 14.58 |
| 200 | 818800 | 46756 | 46552 | 17.59 |
| 300 | 1228200 | 55756 | 55461 | 22.15 |
| 400 | 1637600 | 58048 | 57647 | 28.41 |
| 500 | 2047000 | 59248 | 58750 | 34.84 |
| 700 | 2865800 | 61648 | 60943 | 47.02 |
| 1000 | 4094000 | 65248 | 64230 | 63.74 |
| 1500 | 6141000 | 71248 | 69680 | 88.13 |
| 2000 | 8188000 | 77248 | 75084 | 109.05 |
| ≥314 | 4094·t | 53248+12·t | | |

TABLE VI
OUR METHOD PERFORMANCE FOR GF($2^{14}$)

| t | Chien | mean muls | speedup |
|---|---|---|---|
| 10 | $1.64·10^5$ | $7.64·10^4$ | 2.14 |
| 20 | $3.28·10^5$ | $1.31·10^5$ | 2.50 |
| 30 | $4.91·10^5$ | $1.86·10^5$ | 2.65 |
| 40 | $6.55·10^5$ | $2.40·10^5$ | 2.73 |
| 50 | $8.19·10^5$ | $2.95·10^5$ | 2.78 |
| 70 | $1.15·10^6$ | $4.04·10^5$ | 2.84 |
| 100 | $1.64·10^6$ | $5.68·10^5$ | 2.88 |
| 200 | $3.28·10^6$ | $7.30·10^5$ | 4.49 |
| 300 | $4.91·10^6$ | $7.42·10^5$ | 6.62 |
| 400 | $6.55·10^6$ | $7.55·10^5$ | 8.68 |
| 500 | $8.19·10^6$ | $7.68·10^5$ | 10.67 |
| 700 | $1.15·10^7$ | $7.93·10^5$ | 14.46 |
| 1000 | $1.64·10^7$ | $8.31·10^5$ | 19.72 |
| 1500 | $2.46·10^7$ | $8.94·10^5$ | 27.50 |
| 2000 | $3.28·10^7$ | $9.57·10^5$ | 34.25 |
| 3000 | $4.91·10^7$ | $1.08·10^6$ | 45.40 |
| 5000 | $8.19·10^7$ | $1.33·10^6$ | 61.37 |
| 7000 | $1.15·10^8$ | $1.59·10^6$ | 72.28 |
| ≥128 | 16382·t | 704596+126·t | |

TABLE VII
OUR METHOD PERFORMANCE FOR GF($2^{15}$)

| t | Chien | mean muls | speedup |
|---|---|---|---|
| 10 | $3.28·10^5$ | $2.15·10^5$ | 1.52 |
| 20 | $6.55·10^5$ | $2.62·10^5$ | 2.50 |
| 30 | $9.83·10^5$ | $3.09·10^5$ | 3.18 |
| 40 | $1.31·10^6$ | $3.56·10^5$ | 3.68 |
| 50 | $1.64·10^6$ | $4.03·10^5$ | 4.07 |
| 70 | $2.29·10^6$ | $4.96·10^5$ | 4.62 |
| 100 | $3.28·10^6$ | $6.37·10^5$ | 5.15 |
| 200 | $6.55·10^6$ | $1.10·10^6$ | 5.93 |
| 300 | $9.83·10^6$ | $1.19·10^6$ | 8.25 |
| 400 | $1.31·10^7$ | $1.21·10^6$ | 10.86 |
| 500 | $1.64·10^7$ | $1.22·10^6$ | 13.41 |
| 700 | $2.29·10^7$ | $1.25·10^6$ | 18.32 |
| 1000 | $3.28·10^7$ | $1.30·10^6$ | 25.26 |
| 1500 | $4.91·10^7$ | $1.37·10^6$ | 35.82 |
| 2000 | $6.55·10^7$ | $1.45·10^6$ | 45.29 |
| 3000 | $9.83·10^7$ | $1.60·10^6$ | 61.55 |

| t | Chien | mean muls | speedup |
|---|---|---|---|
| 5000 | $1.64 \cdot 10^8$ | $1.90 \cdot 10^6$ | 86.36 |
| 7000 | $2.29 \cdot 10^8$ | $2.20 \cdot 10^6$ | 104.40 |
| 10000 | $3.28 \cdot 10^8$ | $2.65 \cdot 10^6$ | 123.79 |
| $\geq 216$ | $32766 \cdot t$ | $1146996 + 150 \cdot t$ | |

TABLE VIII
OUR METHOD PERFORMANCE FOR $GF(2^{16})$

| t | Chien | mean muls | speedup |
|---|---|---|---|
| 10 | $6.55 \cdot 10^5$ | $3.06 \cdot 10^5$ | 2.14 |
| 20 | $1.31 \cdot 10^6$ | $4.19 \cdot 10^5$ | 3.13 |
| 30 | $1.97 \cdot 10^6$ | $4.63 \cdot 10^5$ | 4.25 |
| 40 | $2.62 \cdot 10^6$ | $5.07 \cdot 10^5$ | 5.17 |
| 50 | $3.28 \cdot 10^6$ | $5.50 \cdot 10^5$ | 5.95 |
| 70 | $4.59 \cdot 10^6$ | $6.38 \cdot 10^5$ | 7.19 |
| 100 | $6.55 \cdot 10^6$ | $7.69 \cdot 10^5$ | 8.52 |
| 200 | $1.31 \cdot 10^7$ | $1.21 \cdot 10^6$ | 10.87 |
| 300 | $1.97 \cdot 10^7$ | $1.45 \cdot 10^6$ | 13.53 |
| 400 | $2.62 \cdot 10^7$ | $1.48 \cdot 10^6$ | 17.73 |
| 500 | $3.28 \cdot 10^7$ | $1.50 \cdot 10^6$ | 21.78 |
| 700 | $4.59 \cdot 10^7$ | $1.56 \cdot 10^6$ | 29.49 |
| 1000 | $6.55 \cdot 10^7$ | $1.63 \cdot 10^6$ | 40.14 |
| 1500 | $9.83 \cdot 10^7$ | $1.76 \cdot 10^6$ | 55.84 |
| 2000 | $1.31 \cdot 10^8$ | $1.89 \cdot 10^6$ | 69.40 |
| 3000 | $1.97 \cdot 10^8$ | $2.14 \cdot 10^6$ | 91.68 |
| 5000 | $3.28 \cdot 10^8$ | $2.66 \cdot 10^6$ | 123.35 |
| 7000 | $4.59 \cdot 10^8$ | $3.17 \cdot 10^6$ | 144.78 |
| 10000 | $6.55 \cdot 10^8$ | $3.94 \cdot 10^6$ | 166.48 |
| $\geq 254$ | $65534 \cdot t$ | $1376492 + 256 \cdot t$ | |

Starting from $GF(2^{14})$, the difference between the average and maximum number of operations becomes insignificant. Note that this difference is due to optimization (Algorithm 3). It can be seen that the most convenient field for the proposed algorithm is $GF(2^{12})$. Also good results are achieved for $GF(2^{16})$, $GF(2^{15})$, $GF(2^{14})$. In principle, the proposed algorithm has a right to exist for almost all Galois fields, starting from $GF(2^8)$, because of its simplicity.

**REFERENCES**

[1] Chien, R. T., 'Cyclic Decoding Procedures for the Bose-Chaudhuri-Hocquenghem Codes,' IEEE Ref. Trans. on Information Theory, IT-10 (4): 357–363, Oct. 1964, doi:10.1109/TIT.1964.1053699.
[2] E. R. Berlekamp, 'Factoring polynomials over large finite fields,' In SYMSAC '71 - Proceedings of the second ACM symposium on Symbolic and algebraic manipulation, p. 223, New York, USA, 1971. ACM.
[3] S. V. Fedorenko, 'Efficient Algorithm for Finding Roots of Error-Locator Polynomials,' IEEE Access., vol. 9, pp. 38673-38686, Mar. 2021.
[4] T.-K. Truong, J.-H. Jeng, and I. S. Reed, 'Fast algorithm for computing the roots of error locator polynomials up to degree 11 in Reed-Solomon decoders,' IEEE Trans. Commun., vol. 49, no. 5, pp. 779-783, May 2001.
[5] S. V. Fedorenko and P. V. Trifonov, 'Finding roots of polynomials over finite fields,' IEEE Trans. Commun., vol. 50, no. 11, pp. 1709-1711, Nov. 2002.
[6] S. Fedorenko, P. Trifonov, and E. Costa, 'Improved hybrid algorithm for finding roots of error-locator polynomials,' Eur. Trans. Telecommun., vol. 14, no. 5, pp. 411-416, 2003.
[7] T.-C. Lin, T. K. Truong, and P. D. Chen, 'A fast algorithm for the syndrome calculation in algebraic decoding of Reed-Solomon codes,' IEEE Trans. Commun., vol. 55, no. 12, pp. 2240-2244, Dec. 2007.

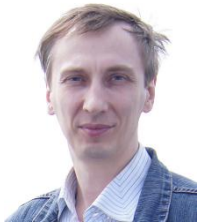

**GENNADY NIKOLAEVICH GLUSHCHENKO** was born in 1969.
He studied at 1986-1995 in MIPT and LPI RAS.
He lives in Moscow.
His main field of activity is microelectronics.